\documentclass[11pt,leqno]{article}
\usepackage{amsfonts, amsthm}
\begin{document}
\newtheorem{exa}{Example}[section]
\newtheorem{theo}{Theorem}[section]
\newtheorem{rem}{Remark}[section]
\newtheorem{lem}{Lemma}[section]
\newtheorem{defi}{Definition}[section]
\newtheorem{coro}{Corollary}[section]
\newtheorem{prob}{Problem}[section]
\def\proj{{\Bbb P}^2} 
\def\R{{\Bbb R}}
\def\Z{{\Bbb Z}}
\def\C{{\Bbb C}}
\def\Q{{\Bbb Q}}
\newcommand{\ffc}[1]{{{\cal F}}(#1)}
\newcommand{\wc}[1]{{{\cal M}}(#1)}
\newcommand{\ps}[1]{{{\Bbb P}}^{#1}}
\def\F{{\cal F}}
\def\P{{\cal P}}
\def\Pt{\C[t]}
\def\RR{{\cal R}} 
\def\li{{\Bbb L}}
\def\L{{\cal L}}
\newcommand{\cfc}[1]{\{\delta_t\}_{t\in #1}}
\begin{center}
{\LARGE\bf Center conditions: Rigidity of logarithmic differential 
equations
\footnote{
\\
Math. classification: 32L30-14D05
\\
Supported by Max Planck Institut F\"ur Mathematik in Bonn}
\\}
\vspace{.25in} {\large {\sc Hossein Movasati}}
\\
{\it Dedicated to Mothers} 
\end{center}
\begin{abstract}
In this paper we prove that any degree $d$ deformation of a generic 
logarithmic polynomial differential equation with a persistent center must be 
logarithmic again. This is a generalization of Ilyashenko's result on
Hamiltonian differential equations. The main tools are Picard-Lefschetz 
theory of a polynomial with complex coefficients in two variables, 
specially the Gusein-Zade/A'Campo's theorem on calculating the Dynkin 
diagram of the polynomial, and the action of Gauss-Manin connection on 
the so called Brieskorn lattice/Petrov module of the polynomial. We will
also
generalize J.P. Francoise recursion formula and (*) condition 
for a polynomial which is a product of lines in a general position.
Some applications on the cyclicity of cycles and the Bautin ideals will
be given. 
 \end{abstract}
\setcounter{section}{-1}
\section{Introduction}
Let us be given the 1-form
\begin{equation}
\label{gravi}
\omega=P(x,y)dy-Q(x,y)dx
\end{equation}
where $P$ and $Q$ are two relatively prime polynomials in 
$\C^2$. The degree of $\omega$ is the maximum of $deg(P)$ and $deg(Q)$. 
The space of degree $d$ $\omega$'s up to multiplication by a 
constant, namely  $\ffc{d}$, is a Zariski open subset of the projective space associated to the coefficient 
space  of polynomial 1-forms ~(\ref{gravi}) with $deg(P),deg(Q)\leq d$.
An element of $\ffc{d}$ induces a holomorphic foliation $\F$ on $\C^2$ 
i.e., the restrictions of $\omega$ to the leaves of $\F$ are 
identically zero. Therefore we denote an element of $\ffc{d}$ by 
$\F(\omega)$ or $\F$ if there is no
confusion about the underlying 1-form in the text and we say that it is of degree $d$.

The points in $Sing(\F(\omega))=\{P=0,Q=0\}$ are called the singularities of $\F$.
 A singularity $p\in\C^2$ of $\F(\omega)$ is called reduced if 
$(P_xQ_y-P_yQ_x)(p)\not =0$. 
A  reduced singularity $p$ is called a center singularity or center for
simplicity if 
there is a holomorphic coordinates system $(\tilde x,\tilde y)$ around $p$ with 
$\tilde x(p)=0,\ \tilde y(p)=0$ and such that in this coordinates system 
$\omega\wedge d(\tilde x^2+\tilde y^2)=0$.
Let $\wc{d}$ be the closure of the subset of $\ffc{d}$ containing $\F(\omega)$'s with at least one center. It is a well-known fact that $\wc{d}$ is
an algebraic subset of $\ffc{d}$ (see for instance \cite{mov3}). 
Now the problem of identifying irreducible components of $\wc{d}$ 
arises.
This problem is also known by the name "Center conditions" in the context of
real polynomial differential equations.

H. Dulac in \cite{dul} proves that $\wc{2}$ has exactly $9$ irreducible 
components (see also \cite{celi} p. 601). 
In this case any foliation in $\wc{2}$ has a Liouvillian first 
integral.
Since this problem finds applications on the number of limit cycles in 
the context of  real differential equations, this classification problem is 
very active. It is recommended to the reader to do a search with the title 
center/centre conditions in mathematical review to obtain many recent 
papers on this problem. We find some partial
results for $d=3$ due to H. Zoladek and others and a similar problem for Abel 
equations $y'=p(x)y^2+q(x)y^3$, $p,q$ polynomials in $x$ (see
also \cite{bfy} and its references).  
In the context of holomorphic foliations we can refer to  
\cite{celi},\cite{muc} and \cite{mov3}. One of the main objectives in this
paper is to introduce some new methods in this problem using an elementary
algebraic geometry. We have borrowed many notions like Brieskorn modules,
Picard-Lefschetz theory and so on from the literature of
singularities of holomorphic functions (see \cite{arn}).   

Let us be given the polynomials $f_i, deg(f_i)=d_i, 1\leq i \leq s$ and 
non-zero complex numbers $\lambda_i\in\C^*, 1\leq i \leq s$. The foliation 
\begin{equation}
\label{apaixonada}
\F=\F(f_1\cdots f_s\sum_{i=1}^s \lambda_i\frac{df_i}{f_i})
\end{equation}
is of degree $d=\sum_{i=1}^{s} d_i-1$ and has the logarithmic first integral 
$f_1^{\lambda_1}\cdots f_s^{\lambda_s}$. 
Let $\L(d_1,\ldots,d_s)$ be the set of foliations ~(\ref{apaixonada}).
Since $\L(d_1,\ldots,d_s)$ is parameterized  by 
$\lambda_i$'s and the coefficients of $f_i$'s, it is irreducible.  
The main theorem of this paper is: 
\begin{theo}
\label{main}
$\overline{\L(d_1,\ldots,d_s)}$ is an irreducible component of $\wc{d}$, where
$d=\sum_{i=1}^{s} d_i-1$ and $\overline{\L(d_1,\ldots,d_s)}$ is the 
closure of $\L(d_1,\ldots,d_s)$ in 
$\F(d)$. 
\end{theo}  
In the case $s=1$ we can assume that $\lambda_1=1$ and so $\L(d+1)$ is 
the space of foliations of the type $\F(df)$, where $f$ is a polynomial of degree $d+1$. This case is proved by Ilyashenko in \cite{ily}. The similar
 result for foliations with a first integral of the type 
$\frac{F^p}{G^q},\frac{deg(F)}{deg(G)}=\frac{q}{p},g.c.d.(p,q)=1$ is 
obtained in \cite{mov}, \cite{mov3}. Some basic tools of this kind of generalizations for Lefschetz pencils on a manifold is worked  in $\cite{muc}$. 

Let us reformulate
our main theorem as follows: Let $\F\in\L(d_1,\ldots,d_s)$ be given by ~(\ref{apaixonada}), $p$
one of the center singularities of $\F$ and $\F_\epsilon$ a holomorphic
deformation of $\F$ in $\ffc{d}$ such that its unique
singularity $p_\epsilon$ near $p$ is still a center. 
There exists an open dense subset $U$ of 
$\L(d_1,\ldots,d_s)$, such that for all $\F\in U$, $\F_\epsilon$ 
admits a logarithmic first integral. More precisely, there
 exist polynomials $f_{i\epsilon},deg(f_{i\epsilon})=d_i,i=1,\ldots,s$
and non-zero complex numbers $\lambda_{i\epsilon}$ such that $\F_\epsilon$ is
given by 
\[
f_{1\epsilon}\cdots f_{s\epsilon} \sum_{i=1}^s \lambda_{i\epsilon}\frac{df_{i\epsilon}}{f_{i\epsilon}}=0
\]
$f_{i\epsilon}$ and $\lambda_{i\epsilon}$ are holomorphic in $\epsilon$ and $f_{i0}=f_i, \lambda_{i0}=\lambda_i, i=1,\ldots,s$.
This new formulation of our main theorem says that the persistence of 
one center implies
the persistence of all other centers and dicritical singularities 
$\{f_i=0\}\cap\{f_j=0\}, i,j=1,\ldots,s$.

We can put $U$ the complement of 
$\L(d_1,\ldots,d_s)\cap sing(\wc{d})$ in $\L(d_1,\ldots,d_s)$. One may not
be satisfied with this $U$ and try to find explicit conditions, for instance:
A foliation $\F(f_1\cdots f_s\sum_{i=1}^s \lambda_i\frac{df_i}{f_i})$ in $U$
satisfies 1. $\{f_i=0\}$ intersects $\{f_j=0\}$ transversally 2. $f_1^{\lambda_1}\cdots f_s^{\lambda_s}$ has nondegenerated critical points in
$\C^2-\cup_{i=1}^s\{f_i=0\}$ and so on. In general one may be interested to 
identify the set $\L(d_1,\ldots,d_s)\cap sing(\wc{d})$.
In any case these questions are not in the focus of this paper.

Since this paper is inspired by Ilyashenko's paper \cite{ily} 
let us give a sketch of the proof in the case $\L(d+1)$:
Let $f$ be a degree $d+1$ polynomial  with the following condition:
1. $f$ has $d^2$ non-degenerate critical points
with distinct values 2. the homogeneous part of $f$ of degree $d+1$ has 
$d+1$ distinct roots.
These conditions are
generic, i.e. in the space of polynomials of degree $d+1$ there is an 
open dense subset such that for all $f$ in this subset the conditions are
satisfied. Let $df+\epsilon\omega+h.o.t.$ be a deformation of $df$ 
such that the singularity  near a center singularity of $df$ , namely 
$p_1$, persists in being center. Then $\int_{\delta}\omega=0$ for all 
vanishing cycles $\delta$ around $p_1$. 
The action of the monodromy on a 
single vanishing cycle $\delta\in H_1(f^{-1}(b),\Z)$, where $b$ is a 
regular value of $f$, generates the whole homology 
(the most significant part of the proof).
Therefore, our integral is zero in all cycles in the fibers of $f$ and 
so it 
is relatively exact.
Knowing the fact $deg(\omega)\leq d$ we conclude that $\omega=dP$, 
where $P$
is a polynomial of degree less than $d+1$. 
Since $\wc{d}$ is an algebraic set, the hypothesis on $df$ is generic
and the tangent vector $\omega$ of any deformation of $df$ in 
$\wc{d}$ is tangent also to ${\cal L}(d+1)$, we conclude that
${\cal L}(d+1)$ is an irreducible component of $\ffc{d}$.
 
Now let us explain our strategy for the proof of Theorem ~\ref{main} and the
structure of the paper. 
First of all, since Picard-Lefschetz theory and classification of
relatively exact 1-forms of 
a multi-valued function $f_1^{\lambda_1}\cdots f_s^{\lambda_s}$ are not well developed, it seems to be difficult to take a generic element of $\L(d_1,\ldots,d_s)$ and then try to repeat Ilyashenko's argument. So 
we look for a special point in $\L(d_1,\ldots,d_s)$. This special point
is going to be $\F_0=\F_0(df),\ f=\Pi_{i=0}^{d}l_i$, 
where $l_i$ is a polynomial 
of degree one in $\R^2$ and the lines $l_i=0$ are in general positions 
in $\R^2$. Every $\L(d_1,\ldots,d_s),\sum_{i=1}^{s}d_i=d+1$ passes 
through $\F_0$ and around $\F_0$ may have many irreducible 
components. 
The main point is to prove that the tangent cone of $\wc{d}$ in $\F_0$ is 
equal to the tangent cone of $\cup_{\sum_{i=1}^{s}d_i=d+1,d_i\in{\Bbb N}\cup\{0\}}\L(d_1,\ldots,d_s)$.
 This will be enough to 
prove our main theorem. To start these calculations three important tools are
needed which I have put them in sections 1,2 and 3. Roughly speaking, in 
\S ~\ref{brieskorn} we want to classify rational 1-forms in $\C^2$ whose
derivatives (Gauss-Manin connection) after some certain order is 
relatively exact. We introduce the
Brieskorn lattice/Petrov module $H$ associated to a polynomial $f$ 
in $\C^2$ and the action of Gauss-Manin connection $\nabla$ on it. 
Using a theorem of 
Mattei-Cerveau
we prove Corollary ~\ref{20mar02} which is enough for our needs in this 
paper. This corollary classifies all $\omega\in H$ with $\nabla^n\omega=0$ for 
a given natural number $n$.
In \S ~\ref{action} we analyze the action  of the monodromy on a Lefschetz 
vanishing
cycle in $f$. 
Using the well-known Theorem ~\ref{bisharm1} and Gusein-Zade/A'Campo's 
Theorem ~\ref{bisharm} we prove Theorem  ~\ref{30mar02} and then Theorem
~\ref{8apr02}.  
In \S ~\ref{melnikov} we consider the deformation $df+\epsilon^k\omega_k+\epsilon^{k+1}\omega_{k+1}+
\cdots+\epsilon^{2k}\omega_{2k} 
+h.o.t.=0,\ \omega_k\not=0, k\in {\Bbb N}$ of $df=0$ with a persistent center. 
We calculate the Melnikov functions $M_i,i=1, \ldots,2k$ and knowing that they
are identically zero we will obtain explicit forms of $\omega_i,i=k,\cdots, 2k$.
\S ~\ref{pufff} is devoted to the proof of Theorem ~\ref{main}. At the end
of this section we discuss our result in the context of real differential equations and
its connection with concepts like cyclicity and Bautin Ideals. 

{\bf Acknowledgments:} Here I would like to thank Max-Planck Institut for Mathematics in
Bonn for financial supporting and for the exceptional ambient of
research. I thank Alcides Lins Neto for putting the paper \cite{ily}
in the first steps of my mathematical carrier. Thanks go also to 
Claus Hertling and  Susumu Tanab\'e for many useful conversations. 
I thank also Konstanz university for hospitality during a visiting.   
\section{Brieskorn lattices/Petrov Modules}
\label{brieskorn}
Let $f$ be a non-composite polynomial of degree $d+1$ in $\C^2$, i.e. $f$ cannot
be composed as $p\circ g$, where $p$ is a polynomial of degree greater than one 
in $\C$ and $g$ is a polynomial in $\C^2$. This condition
is equivalent to the fact that for all $b\in\C$ except a finite number the
fiber $f^{-1}(b)$ is irreducible (see \cite{gom}). Let $\Omega^i,i=0,1,2$ be the set of
polynomial differential $i$-forms in $\C^2$ 
and
$\Pt$ be the ring of polynomials in $t$. $\Omega^i$ is a $\Pt$-module in 
the following way 
\[
p(t).\omega=p(f)\omega,\  p\in\Pt,\omega\in\Omega^i
\]
The Brieskorn lattice/Petrov module  
\[
H=\frac{\Omega^1}{d\Omega^0+\Omega^0 df}
\]
is a $\Pt$-module. Also we define
\[
V=\frac{\Omega^2}{df\wedge\Omega^1}\cong\frac{\C[x,y]}{<f_x,f_y>}
\]
Multiplying by $f$ defines a linear operator on $V$ which is denoted by 
$A$.
\begin{lem}
\label{16apr02}
If $f$ has isolated singularities then the followings are true:
\begin{enumerate}
\item
 $V$ is a $\C$-vector space of 
dimension
$\mu$, where $\mu$ is the sum of local Milnor numbers of $f$;
\item
Eigenvalues of $A$ are exactly the critical values of $f$.
\end{enumerate}
\end{lem}
\begin{proof}
Consider the restriction map $R:V\rightarrow \oplus_p \frac{{\cal O}_{\C^2,p}}{<f_x,f_y>}$, where ${\cal O}_{\C^2,p}$ is the ring of germs of holomorphic functions in a neighborhood of $p$ in $\C^2$ and the sum runs through all critical points  of $f$. 
By Noether's theorem (see \cite{gri} p. 703) $R$ is an isomorphism 
(In fact for surjectivity of $R$ we must modify the proof of Noether's theorem).    Each $\frac{{\cal O}_{\C^2,p}}{<f_x,f_y>}$ is invariant by the linear
operator $A$ and $A-f(p).I$ restricted to it is nilpotent
(see \cite{bri}).
\end{proof}
 
From now on we assume that $f$ has isolated
singularities and we denote the corresponding critical values by 
$c_1,c_2,\ldots,c_r$. 
Let $\tilde{H}$ be the localization of  $H$ by polynomials in $t$ which
vanish only on $c_i$'s and let $p(t)$ be the minimal polynomial 
of $A$. An element of $\tilde{H}$ is a fraction $\omega/a(t),zero(a(t))
\subset
\{c_1,c_2,\ldots,c_r\}$ and we have the usual equality
$\omega/a(t)=\tilde\omega/\tilde a(t)$ if $\tilde a(t)\omega=a(t)\tilde\omega$,
between two fractions.  The
Gauss-Manin connection
\[
\nabla :H\rightarrow \tilde{H}
\]
is defined as follows: 
For an $\omega\in H$ we have $p(f)d\omega=0$ in $V$.
Therefore there is a polynomial 1-form $\eta$ in $\C^2$ such that 
\begin{equation}
\label{18mar02}
p(f)d\omega=df\wedge \eta
\end{equation}
we define $\nabla \omega=\eta/p(t)$. 
Of course we must check that this operator is well-defined. If $\eta_1$ 
and 
$\eta_2$ are two polynomial 1-forms satisfying ~(\ref{18mar02}) then 
$(\eta_1-\eta_2)\wedge df=0$ and so by de Rham lemma 
$\eta_1-\eta_2=Pdf$ ($=0$ 
in $H$),
for a $P$ polynomial in
$\C^2$. Also if $\omega=dP+Qdf$, $P$ and $Q$ two polynomials in $\C^2$, 
then 
$d\omega=dQ\wedge df$ and so $\nabla\omega=dQ$ ($=0$ in $H$).
 
 We can extend $\nabla$ as a function 
from $\tilde{H}$ to $\tilde{H}$ by the rule 
\begin{equation}
\label{10sep02}
\nabla(\omega/q)=(q\nabla\omega-q'\omega)/q^2
\end{equation}
Let $f_1=0,f_2=0,\ldots, f_k=0$ be irreducible components of all 
critical 
fibers of $f$
and $\tilde{\Omega}^i$ be 
the set of rational $i$-forms in $\C^2$ with poles of arbitrary order
along $\{f_i=0\}$'s. 
We define 
$\tilde{H}'=\frac{\tilde{\Omega}^1}{d\tilde{\Omega}^0+\tilde{\Omega}^0 
df}$
and in a similar way as for $H$ a connection $\nabla':\tilde{H}'\rightarrow
\tilde{H}'$ given by the rule ~(\ref{18mar02}).
\begin{lem}
\label{19mar02}
$(\tilde{H},\nabla)$ is isomorphic to $(\tilde{H}',\nabla')$.
\end{lem}
\begin{proof}   
Every rational
1-form in $\C^2$ with poles of arbitrary order
along $\{f_i=0\}$'s determines a unique element of 
$\tilde{H}$ as follows:
if $\omega=\frac{\tilde{\omega}}{f_{i_1}^{r_1}}$ and 
$f-c_1=f_{i_1}^{k_1}\cdots
f_{i_l}^{k_l}$ is the decomposition of $f-c_1$ to irreducible factors 
then 
we multiply both $\tilde{\omega}$ and $f_{i_1}^{r_1}$ by 
$f_{i_1}^{k_1m-r_1}f_2^{k_2m}\cdots f_{i_l}^
{k_lm}$, where $m$ is an integer number satisfying $ 
m-1<\frac{r_1}{k_1}\leq m$, 
and we obtain $\omega=\frac{\tilde{\tilde{\omega}}}{(t-{c_1})^m}$.
Repeating this process by $\tilde{\tilde{\omega}}$ we obtain an element 
of $\tilde{H}$. If $\omega_1=dP+Qdf$, where $P$ and $Q$ are two 
rational
functions on $\C^2$ with poles of arbitrary order
along $\{f_i=0\}$'s  then by applying the above 
method 
on $P,Q$ we can see
that $\omega$ is associated to zero in $\tilde{H}$. Therefore we obtain 
a map
$\tilde{H}'\rightarrow\tilde{H}$ which is the inverse of the canonical 
map
$\tilde{H}\rightarrow\tilde{H}'$ and so it is an isomorphism. Since 
$\nabla$
and $\nabla'$ coincide on $H\subset \tilde{H},\tilde{H}'$, the mentioned isomorphism sends $\nabla'$ to $\nabla$. 
\end{proof}

Let $b$ be a regular value of $f$ and  $\cfc{(\C,b)},\delta_t\in 
f^{-1}(t)$
be a continuous family of cycles in the fibers of $f$.
For an $\omega\in \tilde{H}$ the integral $\int_ \delta\omega$ is 
well-defined and
\begin{equation}
\label{gravdai}
\frac{\partial}{\partial t}\int_{\delta_t}\omega=\int_{\delta_t}\nabla \omega
\end{equation}
(see \cite{arn}). In fact this formula is a bridge between topology and 
algebra in this paper.

Our objective in this section is to analyze the action of $\nabla$ on $H$.
For this purpose let us state a classical theorem. 
Let $\omega$ be a rational 1-form in $\C^2$ and $\cup_{i=1}^k 
\{f_i=0\}$ be
the pole locus of $\omega$.
Suppose that the multiplicity of  $\omega$ along $\{f_i=0\}$ is $r_i$.  
\begin{theo}( \cite{cema})
\label{jahannam}
If $\omega$ is closed, i.e. $d\omega=0$ then there are 
$\lambda_1,\ldots,\lambda_k\in\C$
and a polynomial $g$ such that
\begin{enumerate}
\item
If $r_i=1$ then $\lambda_i\not=0$;
\item
If $r_i>1$ then $f_i$ does not divide $g$;
\item
$\omega$ can be written
\[
\omega=(\sum_{i=1}^{k}\lambda_i\frac{df_i}{f_i})+
d(\frac{g}{ f_1^{r_1-1}\ldots f_k^{r_k-1}})
\]
\end{enumerate}
\end{theo}
Note that if $\omega$ has a pole of order $r_\infty$ at the line at 
infinity
then the degree of $g$ is $\sum_{i=1}^{k}d_j(r_i-1)+r_\infty-1$ and 
$r_\infty+ \sum_{i=1}^{k}\lambda_i d_i=0$.

Let $\L$ be the subset of $\tilde{H}$ containing the 1-forms of the 
type
$\sum_{i=1}^{k}\lambda_i\frac{df_i}{f_i}$, $\lambda_i\in\C$. 
We have the decomposition $\L=\L_1\oplus\cdots\oplus\L_r$, where $\L_i$ 
contains
all the terms $\lambda_j\frac{df_j}{f_j}$, $f_j=0$ being an 
irreducible component of the fiber $f^{-1}(c_i)$. Note that if 
$f^{-1}(c_i)$
is irreducible then $\L_i=0$.  
 
\begin{coro}
\label{20mar02}
For the pair $(\tilde{H},\nabla)$ and a positive integer number $n$
\begin{enumerate}
\item 
\[
Kernel(\nabla^n)=\{\omega\in\tilde{H}\mid \omega=\sum_{j=0}^{n-1}\alpha_jt^j, \alpha_j\in\L\}
\]
\item
\[
Kernel(\nabla^n)\cap H=\{\omega\in H\mid \omega=
\sum_{j=1}^{n-1}\sum_{i=1}^r\alpha_{ij}(t^j-c_i^j),\alpha_{ij}\in \L_i\}
\] 
\end{enumerate}
where $\nabla^n=\nabla\circ\cdots\circ \nabla$ $n$-times.
 \end{coro}
\begin{proof}
Let us prove the first part by induction on $n$. 
We use the isomorphism in Lemma ~\ref{19mar02}. 
If for $\omega\in \tilde{H}$, $\nabla\omega=0$ then $d\omega=dP\wedge 
df$,
where $P$ is a rational function in $\C^2$ with poles along $D$.
Now $d(\omega-Pdf)=0$ and by Theorem ~\ref{jahannam} we have $\omega\in 
\L$. This proves the first part $n=1$ of the induction.
\\
Now if $\nabla^{n+1}\omega=0$ then by induction 
$\nabla\omega=\sum_{j=0}^{n-1}\alpha_jt^j= 
\nabla\sum_{j=0}^{n-1}\frac{\alpha_j}{j+1}t^{j+1}$ or equivalently
$\nabla(\omega-\sum_{j=0}^{n-1}\frac{\alpha_j}{j+1}t^{j+1})=0$. Using the case
$n=1$ we finish the proof of the first part.
\\
Now let us prove the second
part. Let $\omega=\sum_{j=1}^{n-1}\alpha_jt^j+\alpha_0\in H$ and
$\alpha_j=\sum_{i=1}^r\alpha_{ij}$ be the decomposition of $\alpha_j$.
 We write
\[
\omega=\sum_{j=1}^{n-1}\sum_{i=1}^r\alpha_{ij}t^j+\alpha_0=
\]
\[
\sum_{j=1}^{n-1}\sum_{i=1}^r\alpha_{ij}(t^j-c_i^j)+
\sum_{j=1}^{n-1}\sum_{i=1}^r\alpha_{ij}c_i^j+\alpha_0
\]
The first summand in the above belongs to $H$ and hence 
the second one belongs to $\L\cap H$ and so the second summand must be zero.
\end{proof}   
Before we go to the next section let us give three simple but important 
examples. The last one has a very special role in this paper. 
Corollary ~\ref{20mar02} with $n=2$ will be used in the next section. Therefore
we explain it with these examples. 

\begin{exa}\rm
\label{4apr02}
$f=(x^2+y^2-1)x$. Since $xy=0,\ 3x^2x^i=x^i,\  y^2y^i=y^i,i\geq 1$ in 
$V$, $1,x,y,x^2$ form a basis for the vector space $V$. 
$f$ has four critical points
$p_1=(0,1),p_2=(0,-1),p_3=(\sqrt{1/3},0),
p_4=(-\sqrt{1/3},0)$ with three critical values 
$c_1=c_2=0,c_3=-2/3\sqrt{1/3},
c_4=2/3\sqrt{1/3}$.
We define $\omega_1=(x^2+y^2-1)dx$.
We have
\[
 \nabla\omega_1=-\frac{xd(x^2+y^2-1)}{t}=\frac{\omega_1}{t},\
\nabla^2 \omega_1=0
\]
Since $f^{-1}(0)$ is the only reducible fiber of $f$, by Corollary
~\ref{20mar02}, 2 we know that any other 1-form in $H$ with the
property $\nabla^2\omega_1=0$ is some multiply of $\omega_1$ by a 
constant.  
\end{exa}
\begin{exa}\rm
\label{4apr} 
$f=xy(x+y-1)$. There are four critical points 
$p_1=(0,0),p_2=(1,0),p_3=(0,1),
p_4=(1/3,1/3)$ with two critical values $c_1=c_2=c_3=0,c_4=-1/27$. 
Knowing
Corollary ~\ref{20mar02}, we can see that the 1-forms 
$\omega_1=x(x+y-1)dy,
\omega_2=y(x+y-1)dx$ form a basis for the vector space $\{\omega\in 
H\mid
\nabla^2\omega=0\}$.
\end{exa}
\begin{exa}\rm
\label{konstanz}
The lines $l_p=(d-p)x+py-p(d-p)=0, p=0,1,\ldots,d$ are in a general position
in $\R^2$ i.e., they are distinct and no three of them have a common 
intersection point ($l_p$ is the line through $(p,0),(0,d-p)$). 
The polynomial $f=l_0l_1\cdots l_{d}$ satisfies 
the 
properties 

1) All the critical points of $f$ in $\C^2$ are real and non-degenerated;
 
2) The values of $f$ at all saddle critical points equal zero;
\\
By a small perturbation of the lines $l_i$ we also get the property

3) The values of $f$ at center critical points are distinct.
\\ 
(If two critical points associated to two polygons  have the same
value then try to collapse one of the polygons to a point or without volume
region and conclude the above statement. See also Appendix of \cite{mov2} 
for this kind of arguments). 
In a real coordinates system $(\tilde{x},\tilde{y})$ around a saddle
(resp. center) critical point $p$ the function $f$ can be written as 
$f(p)+\tilde x^2-\tilde y^2$ (resp. $f(p)+\tilde x^2+\tilde y^2$).  
$f$ has $a_2=\frac{d(d+1)}{2}$ saddle critical points, $a_1=\sum_{i=2}^{d}[\frac{i-1}{2}]$
center critical points with negative value and 
$a_3=\frac{d(d-1)}{2}-a_1$ center
critical points with positive value, where $[q]$ is the integer number satisfying
$[q]-1<q\leq [q]$. By Corollary ~\ref{20mar02} the set of $\omega\in H$ with
$\nabla^2\omega=0$ is a vector space generated by 
\[
l_0l_1\cdots l_{p-1}l_{p+1}\cdots l_d dl_p, p=0,1,\ldots,d-1
\]
Note that $\sum_{p=0}^d l_0l_1\cdots l_{p-1}l_{p+1}\cdots l_d dl_p=df=0$ in $H$.
\end{exa}
In what follows when we refer to the polynomial $f$ in 
Example ~\ref{konstanz} we mean the one with a small perturbation of the lines
$l_p$ and hence satisfying the property 3.   
\begin{rem}\rm  
 E. Brieskorn in \cite{bri} introduced three 
${\cal O}_{(\C,0)}$-modules $H, H',H''$  in the context of singularity 
of
holomorphic functions $f:(\C^n,0)\rightarrow (\C,0)$ (The $H$ used in 
this paper
is the equivalent of Brieskorn's $H'$). After him these modules
 are called the Brieskorn lattices and recently the similar notions in 
a global
context are introduced by many authors (see \cite{sab},\cite{dim},
\cite{bon},\cite{mov1}).
In the context of differential equations
$H$ appears in the works of G.S. Petrov \cite{pet} on deformation of 
Hamiltonian equations of the type $d(y^2+p(x))$, where $p$ is a 
polynomial in $x$ and
is named by L. Gavrilov in \cite{gav} the Petrov module. 
For this reason I have used both names Brieskorn lattice and Petrov 
module for $H$.

Restriction of an $\omega\in H$ to each fiber defines a global section of
the cohomology fiber bundle of the function $f$ and looking in this way
$\nabla$ is the usual Gauss-Manin connection in the literature. For this
reason I have named $\nabla$ again the Gauss-Manin connection. But of course
we can name $\nabla\omega$ the Gelfand-Leray form of $\omega$ (see \cite{arn}).
\end{rem}     
    
\section{Action of the monodromy}
\label{action}
Suppose that $f$ is a  polynomial function in $\R^2$ with the 
properties 1,2,3 in Example ~\ref{konstanz}. For a $c\in\C$ we define
$L_c=f^{-1}(c)$ in $\C^2$. $\cfc{(\C,c)}$ with $\delta_t\in H_1(L_t,\Z)$ denotes a 
continuous family of cycles.  

Choose a value $b\in\C$ with 
$Im(b)>0$ and  fix a system of paths
joining $b$ with the critical values of $f$, subject to the condition 
that these
paths lie in their entirety in the upper half-plane $Im(z)>0$ except 
for
the ends which coincide with the critical values. Now we can define a
distinguished basis of vanishing cycles in $H_1(L_b,\Z)$ (see \cite{arn}
for the definition).
Critical points of $f$ are in
one to one correspondence with the self intersection points of the real 
curve 
$f=0$, namely $p_j,\ a_1+1\leq j\leq a_1+a_2$, 
and relatively compact components of its complement, namely
$U_i^0, 1\leq i \leq a_1,U_k^2,a_1+a_2+1\leq k\leq a=a_1+a_2+a_3$. 
$U_i^0$ contains a critical point of $f$ with negative 
value and $U_k^2$ contains a critical point with positive value. We
denote by  
\begin{equation}
\label{miqueridatriste}
\delta_i^0, \delta_j^{1},\delta^2_k,\ 1 \leq i \leq a_1,\ 
a_1+1\leq j \leq a_1+a_2,\ a_1+a_2+1\leq k\leq a
\end{equation}
 the distinguished basis of vanishing 
cycles. $\delta^0_i,\delta^2_k$ are called the center vanishing cycles and 
$\delta^1_j$ is called the saddle vanishing cycle. 
\begin{theo}
(S. Gusein-Zade, N. A'Campo)
\label{bisharm}
After choosing a proper orientation for the cycles 
$\delta_i^{\sigma}$  we have
\begin{itemize}
\item
$<\delta_i^{\sigma},\delta_{j}^\sigma>=0$;
\item
$<\delta^0_i,\delta^1_{j}>$ equal to the number of vertices of the 
polygon
$U^0_i$ coinciding with the point $p_j$;  
\item
 $<\delta^2_k,\delta^1_{j}>$ equal to the number of vertices of the 
polygon
$U^2_k$ coinciding with the point $p_j$;
\item
$<\delta^2_k,\delta^0_{i}>$ equal to the number of common edges of $U^2_k$ 
and
$U_i^0$.
\end{itemize}
\end{theo}
This theorem, in an apparently local context, 
is proved by S. Gusein-Zade \cite{gu},\cite{gu1} and N. A'Campo 
\cite{ac},\cite{ac1} independently. However the proof in our case is 
the same. The above theorem gives us the Dynkin diagram of $f$ (see \cite{arn}).

Now let us state another  theorem which we will use in this paper:
\begin{theo}
\label{bisharm1}
In the above situation, the vanishing cycles
$\delta^0_i,\delta^1_j,\delta^2_k$ generate  $H_{1}(L_b,\Z)$ freely. 
\end{theo}
The proof of the above theorems is classical and the reader can consult 
with
\cite{dine},\cite{mov2}. Also the main core of the proof can be found in 
\cite{lam}.

Let us compactify $\C^2$ in $\proj=\{[x;y;z]\mid (x,y,z)\in\C^3-0\}$. 
Here $\C^2=\{[x;y;1]\mid x,y\in \C\}$ and $\li=\{[x;y;0]\mid x,y\in 
\C\}$ 
is the  projective line at infinity. Let $f=f_0+f_1+\cdots+f_{d+1}$ be 
the 
decomposition of $f$ to homogeneous parts.  We look at $f$ as a 
rational function on $\proj$ by rewriting $f$ as
\[
f=f(x/z,y/z)=\frac{z^{d+1}f_0(x,y)+z^df_1(x,y)+\cdots+f_{d+1}(x,y)}{z^{d+1}}
\]
The indeterminacy set of $f$ is given by 
\[
\RR=\{[x;y;0]\mid f_{d+1}(x,y)=0\}
\]
Now suppose that $\RR$ has $d+1$ distinct points. For instance the polynomial in
Example ~\ref{konstanz} has this property.
 This implies
that the fibers of $f$ intersect the line at infinity transversally.
Doing just one blow-up in each point of $\RR$ and using Ehresmann's 
fibration theorem we conclude that the
map $f$ is a $C^{\infty}$ fibration on $\C-C$, where $C=\{c_1,\ldots,c_r\}$ is the set of critical values of $f$. 
In general case we must do more
blow-ups to obtain this conclusion and the set of critical points of $f$ may be
a proper subset of $C$. Therefore we have the action of the
monodromy on the first (co)homology group of $L_b$:
\[
h:\pi_1(\C-C,b)\times H_1(L_b,\Z)\rightarrow H_1(L_b,\Z)
\]
Recall the  system of paths in the beginning of this section. When we say the
"Monodromy around a critical value " we mean the monodromy associated to the
path which gets out of $b$, goes along $\gamma$ (the path connecting $b$ to the
critical value in this system of paths), turns around the critical value 
counterclockwise and then comes back to $b$ along $\gamma$.
By Picard-Lefschetz formula (see \cite{lam}) the monodromy around zero 
is given by
\begin{equation}
\label{azizam}
\delta\rightarrow \delta-\sum_{j=a_1+1}^{a_2}<\delta,\delta^1_j>\delta^1_j
\end{equation}
and the monodromy around a center critical value $c_i$ is given by
\begin{equation}
\label{azizam1}
\delta\rightarrow \delta-<\delta,\delta_i>\delta_i
\end{equation}
For the regular value 
$b\in\C-C$, $\overline{L_b}=L_b\cup \RR$ is a compact Riemann
surface. Let $I$ be the subgroup of $H_1(L_b, \Z)$ generated by the 
cycles
around the points of $\RR$ in $\overline{L_b}$. We have
\begin{equation}
\label{azizam2}
I=\{\delta\in H_1(L_b, \Z)\mid <\delta,\delta'>=0 ,\ \forall 
\delta'\in
H_1(L_b, \Z) \}
\end{equation}
Elements of $I$ are fixed under the action of the monodromy. 

\def\cc{{\cal C}}
Now we are in a position such that we can look at our objects in an abstract way:
We have a union of curves $\cc=f^{-1}(0)$ in $\R^2$. To $\cc$  we associate an Abelian group
$G$ ($=H_{1}(L_b,\Z)$) freely generated by the 
symbols ~(\ref{miqueridatriste}). These symbols are in one to one correspondence with the self intersection points of $\cc$ and the relatively compact components of the complement of $\cc$ in $\R^2$. We have an 
antisymmetric pairing $<.>$ given by the items of  
Theorem ~\ref{bisharm}. Also the non-Abelian freely generated group 
$\pi_1(\C-C,b)$ acts from left on $G$ with the rules
$~(\ref{azizam}),~(\ref{azizam1})$. We have also the subgroup $I$ of $G$ defined by $~(\ref{azizam2})$. These are all we are going to need. From now on we can
think about (vanishing) cycles in this abstract context.   

We will apply the above arguments for the Example ~\ref{konstanz}. 
For the  line $l_p=0, p=0,1,\ldots, d$ we can associate 
the saddle critical points of $f$ on $l_p=0$ and the corresponding
vanishing cycles. We rename these vanishing cycles
by  $\delta_{j}^{l_p}, j=1,2,\ldots, d$ and suppose that the ordering 
by 
the index $j$ is the same as the ordering of corresponding saddle 
points
in the line $l_p=0$ 
(there are two ways for such indexing, we choose one of them). We define
\[
\delta^{l_p}=\sum^d_{j=1}(-1)^j \delta_{j}^{l_p}\in G, p=0,1,\ldots, 
d
\]
Now reindex all relatively compact polygons by $U_i, 1\leq i\leq 
a_1+a_3$. 
For any polygon $U_i$ we denote by $\delta^{i}$ ($\in G$) the sum of 
vanishing cycles in the vertices of $U_i$. Also we denote by $\delta_i$ 
the vanishing cycle associated to $U_i$ 
\begin{lem}
\label{gonah}
The cycles  $\delta^{l_p},\ 1\leq p \leq d$ and $\delta^{i},1\leq i 
\leq a_1+a_3$ generate
all saddle vanishing cycles in $G\otimes \Q$ and so
they are linearly independent in $G\otimes \Q$. 
\end{lem}
Note that the number of the cycles considered in the lemma is equal to the
number of saddle vanishing cycles.    
\begin{proof}
This is a nice high school problem. It would be more difficult if we assume 
only that the lines $l_p$ are in a general position
in $\R^2$, i.e. no three of them have a common 
intersection point and no two of them are parallel. 
For our example we give the following hint:
1. First draw  the lines for a small value of $d$ 2. Let 
$\delta_{p, p+1}$ denote the vanishing cycle associated to the intersection of $l_p$ and $l_{p+1}$, $p=0,1,\ldots, d-1$. Try to write
$d. \delta_{p, p+1}$ as 
an integral sum of $\delta^{l_p}$ and $\delta^{l_{p+1}}$ and $\delta^i$, 
where $i$ runs through the index of all polygons between (the angle less than $90^\circ$) the lines $l_p$ and $l_{p+1}$. 3. Now it is easy to conclude that
every vanishing cycle associated to the intersection points multiplied by $d$
can be written as an integral sum of $\delta^{l_p},\ 0\leq p \leq d$ and $\delta^{i},1\leq i 
\leq a_1+a_3$. 4. After choosing a proper sign for 
$\delta^{l_p},\ 0\leq p \leq d$ prove that $\sum_{i=0}^d\delta^{l_p}=0$. 
 \end{proof}  

Now let us state the geometric meaning of $\delta^i$ and $\delta^{l_p}$.
 \begin{lem}
We have
\begin{enumerate}
\item
$\delta^{i}= \delta_i-h_0(\delta_i)$, where $h_0$ is the monodromy around 
$0$;
\item
 $I\otimes \Q$ is generated by the cycles $\delta^{l_p}, 1\leq p \leq d$.
\end{enumerate}
\end{lem}
\begin{proof}
The first part is a direct consequence of Picard-Lefschetz formula and 
Theorem
~\ref{bisharm}.
\\
By Theorems ~\ref{bisharm} and ~\ref{bisharm1} we have 
$<\delta^{l_p},\delta'>=0 ,\ \forall \delta'\in
H_1(L_b, \Z)$ and so $\delta^{l_p}$ is in $I$. 
By Lemma ~\ref{gonah} $\delta^{l_p},k=1,2,\ldots, d$ are linearly 
independent in $G$ and we know that $I$ is freely generated of rank 
$d$.
Therefore $\delta^{l_p},k=1,2,\ldots, d$ generate  $I\otimes \Q$ freely.   
\end{proof}
Because of the symmetry in Example ~\ref{konstanz}, one may conjecture
that $\delta^{l_p}$ is the cycle around $\{l_p=0\}\cap\li$ in 
$\overline{L_b}$ (multiplied by a rational number). Since we don't need this statement we don't try to prove it.

\begin{theo}
\label{30mar02}
In the Example ~\ref{konstanz} the action of the monodromy on a Lefschetz 
 vanishing cycle generates the homology 
$H_1(\overline{L_b},\Q)$.
\end{theo}
\begin{proof}
First consider the case where $\delta$ is a center vanishing cycle. 
By Theorem ~\ref{bisharm} and Picard-Lefschetz formula the action of the
monodromy generates $\delta_i$ and then $\delta_i$'s. By Lemma ~\ref{gonah}
and Theorem ~\ref{bisharm1} these cycles generate 
$H_1(\overline{L_b},\Q)$. 
Now suppose that $\delta$ is a saddle vanishing cycle.
There is a center vanishing cycle $\delta'$ such that $<\delta',\delta>\not=0$.
Performing a monodromy of $\delta$ around the critical value associated to $\delta'$ and subtracting the obtained cycle by $\delta$ we obtain $\delta'$ and
so we fall in the first case.
\end{proof}

In the beginning of this section we defined the degree of a polynomial 1-form
$\omega=Pdy-Qdx$ to be the maximum of $deg(F)$ and $deg(Q)$. This definition
is no more useful when we look at $\omega$ as a rational 1-form in $\proj$
or when we consider the foliation induced by $\omega$ in $\proj$. Let us
 introduce a new definition of degree.  For a polynomial 1-form $\omega$ we define $deg_1(\omega)$ to be the 
order of
$\omega$ along the line at infinity mines two. We can see easily that 
if 
$deg_1(\omega)\leq d$ then $\omega=Pdy-Qdx+ G(xdy-ydx)$, where $P,Q$ 
are two 
polynomials of degree at most $d$ and $G$ is zero or a homogeneous 
polynomial
of degree $d$. Therefore $deg_1(\omega)-deg(\omega)=0,1$.  
Naturally for a  $\omega\in H$ we define $deg_1(\omega)$ to be the 
minimum of the $deg_1$'s of the elements of $\omega$. Let $q$ be an indeterminacy
 point of $f$ at the line at infinity $\li$. Recall that the fibers of $f$ intersect $\li$ transversally. Now we can choose a continuous family
of cycles $\cfc{\C}$ such that $\delta_t$ is a cycle in $L_t$ around $q$. 
Latter we will need
the following lemma.
\begin{lem}
\label{16jan02}
For $\omega\in H$ the integral $\int_{\delta_t}\omega$ as a function in $t$ 
is a polynomial
of degree at most $[\frac{n}{d+1}]$, where $n-2=deg_1(\omega)$ and 
$d+1=deg(f)$.
$\nabla^i\omega,
i>[\frac{n}{d+1}]$ restricted to each fiber $\overline{f^{-1}(t)}$ has 
not residues in $\RR$ and 
hence is a 1-form of the second type.
\end{lem}
\begin{proof}
We have $p(t):=\int_{\delta_t}\omega=
t^\frac{n}{d+1}\int_{\delta_t}\frac{\omega}{f^{\frac{n}{d+1}}}$. Since 
the 
1-form 
$\frac{\omega}{f^{\frac{n}{d+1}}}$ has not pole along the line at infinity,
 $\frac{p(t)}{t^\frac{n}{m}}$ has finite growth at $t=\infty$. Since 
$p(t)$ is holomorphic in $\C$, we conclude 
that $p(t)$ is a polynomial of degree at most
$[\frac{n}{d+1}]$. 
The second part is a direct consequence of the first one and the 
formula ~(\ref{gravdai}).
\end{proof}
Let $f$ be a polynomial and $\omega$ be a 1-form in $\C^2$. $\omega$ is called
relatively exact modulo $f$, or simply relatively exact if the underlying $f$
is known, if the restriction of $\omega$ to each fiber $L_b$ is exact.
1-forms of the type $dP+Qdf$, where $P$,$Q$ are polynomials in $\C^2$,
are relatively exact. Latter we will need the following lemma:
\begin{lem}
\label{6apr02}
If the fibers $f^{-1}(t),t\in\C$ of a polynomial $f$ are topologically
connected then every relatively exact 
1-form $\omega$ is of the form
$dP+Qdf$, where $P,Q$ are two polynomials in $\C^2$ with
$deg(P)=deg_1(\omega)+2$ and $\deg(Q)=deg_1(\omega)+2-deg(f)$.
\end{lem} 
The proof can be found in \cite{gav, bon, bon1}. 
For the assertion about the degrees see \cite{mov3} Theorem 4.1.
This kind of results was obtained for the first time by Ilyashenko
in \cite{ily}. 
The main objective of this section is the following theorem.
\begin{theo}
\label{8apr02}
In the Example ~\ref{konstanz} let $\delta_t$ be a continuous family of vanishing cycles and $\omega$ be a degree $d$ 1-form in $\C^2$ 
such that $\int_{\delta_t}\omega=0,t\in(\C,b)$. Then $\omega$ is of the form
\[
\omega=l_0\cdots l_d \alpha+d(P),\ 
\alpha=\sum_{i=0}^{d}\lambda_i\frac{dl_i}{l_i}, \lambda_i\in\C
\]
where $P$ is a polynomial of degree not greater than $d+1$.
\end{theo}
The statement of the above theorem for $f$ can be considered as $(*)$ condition 
of J.P. Francoise in \cite{fra}.
\begin{proof}
By the hypothesis and ~(\ref{gravdai}) we have 
$\int_{\delta_t}\nabla^2\omega=0,  t\in(\C,b)$.
Lemma ~\ref{16jan02} implies that  $\nabla^2\omega$ is a 1-form without 
residue in each fiber and Theorem ~\ref{30mar02} implies 
that $\nabla^2\omega$ is a relatively exact 1-form and hence by
Lemma ~\ref{6apr02} it is zero in $\tilde{H}$ (Note that by ~(\ref{18mar02}) 
and  ~(\ref{10sep02}) $\nabla^2\omega$ is
of the form $\eta/p(t)^2,\eta\in H$ and hence $\eta$ is relatively exact). 
Since $L_0$ is the only reducible fiber of $f$, 
by Corollary ~\ref{20mar02} $\omega$ must be of the form
$f\sum_{i=0}^{d}\lambda_i\frac{dl_i}{l_i}+dP+Qdf$, where $P$ and $Q$ are two
polynomials in $\C^2$ and $\lambda_i$ is a complex number. Recall that
$\omega$ has degree  $d$.
By Lemma ~\ref{6apr02} we can suppose that $P$ (resp. $Q$) has degree less 
than $d+2$ (resp. $1$).  
We have $dP_{d+2}+Q_{1}df_{d+1}=0$, where $P_{d+2}$ denote the
homogeneous part of $P$ of degree $d+2$ and so on. If
$Q_{1}$ is not identically zero then this equality implies that 
$f_{d+1}={Q_1}^{d+1}$ which is not 
our case (write $d(P_{d+2}+Q_1f_{d+1})-f_{d+1}dQ_1=0$ and then conclude that
$P_{d+2}$ and $f_{d+1}$ are polynomials in $Q_1$. Since $f_{d+1}$ is homogeneous in $x$ and $y$ it must be of the claimed form). 
Therefore $Q$ must be 
constant and $P$ of degree less than $d+1$. We substitute $P+Qf$ by $P$
and so we can assume that $Q$ is zero. We obtain the desired equality.
\end{proof}

\section{Deformation and Melnikov functions}
\label{melnikov}
As a first attempt to prove Theorem ~\ref{main} one may fix a generic  
$\F\in {\cal L}(d_1,\ldots,d_s)$ and perform a deformation 
such that one of the center singularities of $\F$ persists. 
Then one may try to find some tools for finding a logarithmic first integral 
for the deformed foliation. 
These tools in the Ilyashenko's case ${\cal L}(d+1)$ were Picard-Lefschetz
 theory of a polynomial and the classification of relatively exact 1-forms modulo the polynomial. Developing these tools for a generic point of other irreducible components of $\wc{d}$, for instance ${\cal L}(d_1,\ldots,d_s)$, seems to be 
difficult and wasting the time. 

In this section we want to explain this idea that it is not necessary 
to take a generic point of ${\cal L}(d_1,\ldots,d_s)$.
For instance suppose that the variety 
${\cal L}(d_1,\ldots,d_s)$ has a point $\F(df)$ which is Hamiltonian. 
This point may lie in other irreducible components of $\wc{d}$. 
 Now the idea is to deform $\F(df)$ in such a way that one of its 
centers persists and since we can develop our tools for the Hamiltonian system 
$\F(df)$, we can calculate the tangent cone of $\wc{d}$ in $\F(df)$.
Now if the tangent cone of one of the branches of ${\cal L}(d_1,\ldots,d_s)$ at
$\F(df)$ is an irreducible component of the tangent cone of $\wc{d}$ at
$\F(df)$ then
${\cal L}(d_1,\ldots,d_s)$ is  irreducible component of $\wc{d}$ locally and
hence globally.
 
Now we are going to realize this idea for the Hamiltonian foliation
with the polynomial in Example ~\ref{konstanz}. 
This will lead to the proof of our main theorem. 
In this section  $\L$ denotes the set of rational 1-forms of the type 
$\sum_{i=0}^{d}\lambda_i\frac{dl_i}{l_i}$ and $\P_n$ (resp. $\P_*$) denotes the set of
polynomials of degree not greater than $n$ (resp. arbitrary degree) in $\C^2$.
 
Let $f$ be the polynomial considered in Example ~\ref{konstanz} and
\begin{equation}
\label{taadol}
\F_\epsilon: \omega_\epsilon=df+\epsilon^k\omega_k+\epsilon^{k+1}\omega_{k+1}+
\cdots+\epsilon^{2k}\omega_{2k} 
+h.o.t.=0,\ \omega_k\not=0
\end{equation}
where $\omega_i, i=k,k+1,\ldots$ are polynomial 1-forms in $\C^2$ and 
$k$ is a natural number, 
be a one parameter degree $d$ deformation of $\F(df)$.
\begin{defi}
Let $p$ be a center 
singularity of $\F=\F(df)$. $p$ is called a persistent center if there exists
a sequence $\epsilon_i,i=1,2,3,\ldots,\ \epsilon_i\rightarrow 0$ such that the singularity of $\F_{\epsilon_i}$ near $p_0=p$, namely $p_{\epsilon_i}$,
is center for all $i$.
\end{defi}
Latter we will see that if $p$ is persistent then $p_\epsilon$ is center
for all $\epsilon$.
Let $I=\{0,1,\ldots,d\}$. For an equivalence relation $J$ in $I$ we denote by
 $J_1,J_2,\ldots, J_{s_J}$ all equivalence classes of $J$ and we define
$f^J_i=\Pi_{j\in J_i}l_j,
i=1,2,\ldots,s_J$. Our main theorem in this section is the following: 
\begin{theo}
\label{mude}
If $p$ is a persistent center in the degree $d$ deformation 
~(\ref{taadol}) then $\omega_k$ is
of the form
\[
\omega_k=l_0l_1\cdots l_d\sum_J \sum_{i=1}^{s_J}
(\lambda^J_i\frac{df^J_i}{f^J_i}+\frac{A^J_i}{f^J_i})
\]
where  in the first sum $J$ runs through all equivalence relations in $I$,
for each $J$ the complex numbers $\lambda^J_i,i=1,2,\ldots,s_J$ are distinct
and $A^J_i\in\P_{deg(f^J_i)}$.
\end{theo}

Let  $\delta_t,t\in(\C,b)$ be a continuous family of 
vanishing cycles around $p$ and $\Sigma$ be a transverse section to 
$\F$ at 
some point of $\delta_b$. We write the Taylor expansion of the 
deformed holonomy $h_\epsilon(t)$
\[
h_\epsilon(t)-t=M_1(t)\epsilon+M_2(t)\epsilon^2+\cdots+M_i(t)\epsilon^i+h.o.t.
\]
$M_i(t)$ is called the $i$-th Melnikov function of the deformation.
If the center $p$ is persistent under the deformation then $M_i= 0$ for
all $i$. But we don't need to use all these equalities. For instance in
Ilyashenko's case ${\cal L}(d+1), k=1$ we need only $M_1=0$. 
To prove Theorem ~\ref{mude} and our main theorem we will need to use $M_k=M_{k+1}=\cdots=M_{2k}=0$.
\begin{lem}
Let $M_i$ be the $i$-th Melnikov function 
associated to the deformation ~(\ref{taadol}). We have: 
1) $M_1=M_2=\cdots=M_{k-1}= 0$ 2) if $\Sigma$ is parameterized 
by the image of $f$, i.e. $t=f(z),z\in\Sigma$ then 
\[
M_k(t)=-\int_{\delta_t}\omega_k
\]
3) If $M_k=M_{k+1}=\cdots=M_{2k-1}
= 0$ then 
\[
\omega_i=f\alpha_i+dP_i,\ P_i\in\P_{d+1},\ \alpha_i\in\L, \  k\leq i \leq 2k-1
\]
4) Moreover if the transverse
 section is parameterized by the image of a branch of $lnf$, i.e. 
$t=lnf(z),z\in\Sigma$ then
\[
M_{2k}(t)=-\int_{\delta_t} \frac{\omega_{2k}}{f}- \frac{P_k}{f}.\alpha_k
\]   
\end{lem}
This theorem can be considered as a generalization of Francoise recursion
formula (see \cite{fra}). Note that for our polynomial, which is a product of lines, we
have Theorem ~\ref{8apr02} instead of Francoise (*) condition.
\begin{proof}
Let $\delta_{t,h_{\epsilon}(t)}$ be the path connecting $t$ and 
$h_{\epsilon}(t)$
along $\delta_b$ in the leaf of $\F_\epsilon$ through $t$. We take
the integral $\int_{\delta_{t,h_{\epsilon}(t)}}$ of ~(\ref{taadol}). Now the
equalities associated to the coefficients of $\epsilon^i, 1\leq i\leq k-1$ prove
the first part. The equality associated to the coefficient of $\epsilon^k$ proves
the second part (for more detail see \cite{rou}\cite{mov}).
\\
We prove the third and fourth part by induction on $i$. First $i=k$. 
$M_k=0$ implies that $\int_{\delta_t}\omega_k=0,\ t\in(\C,b)$ and so
by Theorem ~\ref{8apr02} $\omega_k$ is of the form   
\begin{equation}
\omega_k=f\alpha_k+dP_k,\ \alpha_k\in\L,
P_k\in\P_{d+1}
\end{equation}
Now let us suppose that
\[
\omega_j=f\alpha_j+dP_j,\ \alpha_j=\sum_{p=0}^{d}\lambda_{j,p}\frac{dl_p}{l_p}
,\ k\leq j\leq i
\]
Let $\overline\omega_\epsilon=\omega_\epsilon/f$, 
$\overline{\omega}=\omega/f$ and
so on. With this new notation we have
\[
\tilde\omega_j=\alpha_j+ \overline {dP_j}=d(\overline P_j+ ln l_1^{\lambda_{j,1}}\cdots 
l_{d+1}^{\lambda_
{j,d+1}})+\overline P_j \overline{df}=dg_j+\overline P_j \overline {df}
\]
From now on suppose that $\Sigma$ is parameterized by the image of $lnf$.
We have
\[
(1-\overline P_k\epsilon^k)\cdots(1-\overline P_i\epsilon^i)
\overline\omega_\epsilon
\]
\[
(1-\overline P_k\epsilon^k)\cdots(1-\overline P_i\epsilon^i)
(\overline{df}+\epsilon^k\overline \omega_k+\epsilon^{k+1}\overline\omega_{k+1}+
\cdots+\epsilon^{2k}\overline \omega_{2k} 
+h.o.t.)
\]
\begin{equation}
\label{9apr02}
=\overline {df}+\epsilon^k dg_k+\epsilon^{k+1}dg_{k+1}+\cdots+
\epsilon^i dg_i 
\]
\[
+\epsilon^{i+1}\overline\omega_{i+1}+
\cdots+\epsilon^{2k-1}\overline \omega_{2k-1}+ 
\epsilon^{2k}(\overline \omega_{2k}-\overline P_k\overline \omega_k)
+h.o.t.
\end{equation}
 We take
the integral $\int_{\delta_{t,h_{\epsilon}(t)}}$ of ~(\ref{9apr02}). 
Now the
equality associated to the coefficients of $\epsilon^{i+1}$ is 
$M_{i+1}(t)+\int_{\delta_t}\overline\omega_{i+1}=0$. $M_{i+1}= 0$ and Theorem
~\ref{8apr02} imply that $\omega_{i+1}$ is of the desired form. 
In the last step $i=2k-1$ the fourth part of the lemma is proved.
Note that $\int_{\delta_t}\overline \omega_{2k}-\overline P_k\overline \omega_k=
\int_{\delta_t}
\overline \omega_{2k}-\overline P_k\alpha_k$.
\end{proof}

Now $M_{2k}=0$ implies that $\int_{\delta_t}f\omega_{2k}-P_kf\alpha_k=0$.
$\omega_{2k}-P_kf\alpha_k$ has a pole of order at most $2d+3$ at the line at infinity. Therefore by Lemma   ~\ref{16jan02} $\nabla^3(\omega_{2k}-
P_kf\alpha_k)$ restricted to the fibers $\overline{f^{-1}(b)}$ has not residues.
By Theorem ~\ref{30mar02} and Corollary ~\ref{20mar02} we conclude that
\[
f\omega_{2k}-P_k . f\alpha_k=\alpha_1f +\alpha_2f^2+dg+pdf, \ g,p\in\P_*,\  \alpha_i\in\L,\ i=1,2
\]
The restriction of the above equality to the $L_0=f^{-1}(0)$ implies that
$g$ is constant on  $L_0$. Since $L_0$ is connected in $\C^2$, we conclude that 
$dg=d(fg'),g'\in\P_*$. From now on we write $\lambda_{k,i}=\lambda_i$. 
The above equality modulo multiplications by $l_i$ gives us
\begin{equation}
\label{khaste}
l_i\mid \lambda_{i} P_k+g',\ \ i=0,1,\ldots, d
\end{equation}
Let $I=\{0,1,\ldots,d\}$. Define  $i\cong j$ if $\lambda_i=\lambda_j$. 
$\cong$
is an equivalence relation. Let $J_1,J_2,\ldots,J_s$ be the equivalence
classes of $\cong$. We define $f_i=\Pi_{j\in J_i}l_j$. Note that
we have $f=f_1f_2\cdots f_s$. The following lemma finishes the proof of
Theorem ~\ref{mude}. 
\begin{lem}
In the above situation, $P_k$ must be of the form
\begin{equation}
\label{19apr02}
P_k=f\sum_{i=1}^s \frac{A_i}{f_i},\  A_i\in\P_{\# J_i}
\end{equation}
\end{lem}
\begin{proof}
Let $d_i=deg(f_i)$. ~(\ref{khaste}) implies that $P_k$ is zero in 
$\{f_i=0\}\cap \{f_j=0\}$. The space of $P\in\P_{d+1}$ vanishing in $\{f_i=0\}\cap \{f_j=0\}$ for all $1\leq i<j\leq s$, namely $G$, 
is of dimension $\frac{(d+2)(d+3)}{2}-\sum_{1\leq i<j\leq s} d_id_j$. 
(The matrix $[P_m(B_n)]$ where $P_m$ runs in $\P_{d-1}$ and $B_m$ in the intersection points of the lines $l_i$ has non zero determinant, otherwise there would be a polynomial $P$ of degree not greater than $d-1$ vanishing 
in all $B_m$ which is a contradiction, because $P=0$ intersects a line at most
in $d-1$ points. Therefore the map 
$\psi: \P_{d+1}\rightarrow \C^{d(d+1)/2},\  
\psi(P)=(P(B_m))$ is surjective and hence the map $\psi':\P_{d+1}\rightarrow \C^{\sum_{1\leq i<j\leq s} d_id_j}, 
\psi'(P)=(P(B_m))$ is surjective, where in the second map $B_m$ runs in
the intersection points of $f_i$'s). But the space of
polynomials in ~(\ref{19apr02}) is a subset of $G$ and has dimension
$-1+\sum_{i=1}^s\frac{(d_i+1)(d_i+2)}{2}$. Since $d+1=\sum_{i=1}^{s}d_i$, these 
two numbers are equal  and so
$P_k$ must be of the form ~(\ref{19apr02}).
\end{proof}

\section{Proof of Theorem ~\ref{main} }
\label{pufff}
\newcommand{\tc}[2]{TC_{#1}{#2}}
Let $(X,0)$ be a germ of an analytic variety in $(\C^n,0)$. The tangent
cone $\tc{0}{X}$ of $X$ at $0$ is defined as follows:
Let $\gamma:(\C,0)\rightarrow (\C^n,0)$ be an analytic map such that 
its image lies in $X$ and has the Taylor series $\gamma=
\omega\epsilon^l+\omega'\epsilon^{l+1}+\cdots,\omega,\omega',\ldots\in\C^n$. 
$T_l$ is the set of all $\omega$ and $\tc{0}{X}=\cup_{i=1}^\infty T_l$.

We have $\C .\tc{0}{X}=TC_0X$ therefore we can projectivize $\tc{0}{X}$ 
and obtain a subset, namely $Y$, of $\ps{n-1}$.
Suppose that $X$ is irreducible. 
Let $\pi: M\rightarrow (\C^n,0)$ be a blow-up at $0$ with the divisor 
$\pi^{-1}(0)\cong \ps{n-1}$. The closure  $\tilde{X}$ of $\pi^{-1}(X-\{0\})$ 
in $M$ is an irreducible analytic set and we can see easily that 
$Y$ is isomorphic to the intersection of  $\tilde{X}$ and $\ps{n-1}\subset M$,
and so it is algebraic compact subset of $\ps{n-1}$.
Moreover since $dim(Y)\geq dim(\tilde{X})+dim\ps{n-1}-n$ 
(see \cite{kem} Theorem 3.6.1) and $Y$ cannot be the whole $\ps{n-1}$, 
$Y$ is of  pure dimension $dimX-1$, i.e. each
irreducible component of $Y$ is of dimension $dimX-1$. We conclude that
$\tc{0}X$ is an algebraic subset of $\C^n$ of pure dimension $dim(X)$.
If $X$ is smooth then $\tc{0}{X}$ is the usual tangent space of $X$ at $0$ and
hence it is a vector space. For more information about the tangent cone
of a singularity and its definition  by the leading terms of the polynomial
defining the singularity see \cite{mum} and \cite{kem} Section 6.2.

The variety $\L(d_1,\cdots,d_s)$ is parameterized by 
\[
\tau:\C^s\times \P_{d_1}\times\cdots\times\P_{d_s}\rightarrow\ffc{d}
\]
\begin{equation}
\label{29apr02}
\tau(\lambda_1,\ldots,\lambda_s,f_1,\ldots,f_s)= 
f_1\cdots f_s\sum_{i=1}^s \lambda_i\frac{df_i}{f_i}
\end{equation}
and so it is irreducible. Let $J$ be an equivalence relation in
$I=\{0,1,\ldots, d\}$ with $s$ equivalence classes, namely $J_1,\ldots,J_s$.
Let also $f$ be the polynomial in Example ~\ref{konstanz} and $\F_0=\F(df)$.
In a neighborhood of $\F_0$ in $\ffc{d}$, $\L(d_1,\cdots,d_s)$ has many irreducible 
components corresponding 
to the $J$'s as follows: 

The above parameterization  near 
$(1,\ldots,1,\Pi_{i\in J_1}l_i,\ldots,\Pi_{i\in J_s}l_i)$ determines an irreducible
component, namely $\L(d_1,\cdots,d_s)_J$, of $(\L(d_1,\cdots,d_s),\F_0)$ corresponding to $J$. Fix one of these branches and name it $X$.  
Now to prove our main theorem  it is enough to prove that $X$ is an 
irreducible component of $(\wc{d},\F_0)$.    

What we have proved in Theorem ~\ref{mude} is:
\[
\tc{\F_0}{\wc{d}}=\cup \tc{\F_0}{\L(d_1,\cdots,d_s)_J}
\]
where the union is taken over all $1\leq d_i\leq d+1,
\ 1\leq i\leq s\leq d+1,\ \sum_{i=1}^{s}d_i=d+1$ and all equivalence relations
$J$. Now let  $X\subset X'$, where $X'$ 
is an irreducible component of $(\wc{d},\F_0)$.
Since the above union is the decomposition of $\tc{\F_0}{\wc{d}}$ to 
irreducible components, the  irreducible 
component of $\tc{\F_0}{X'}$ containing $\tc{\F_0}{X}$ must be a subset of
$\tc{\F_0}{X}$ and so is equal $\tc{\F_0}{X}$. An irreducible component of
$\tc{\F_0}{X'}$ is of dimension $dim(X)$ and so $X'$ is of dimension $dim(X)$.
Since $X\subset X'$ and $X,X'$ are irreducible, we conclude that $X=X'$.\qed

{\bf Limit cycles and Bautin Ideals:}
Let $\F_0\in\wc{d}$, 
$p_0$ be a center singularity of $\F_0$, $\delta_t,t\in(\C,0)$ be a 
continuous family of cycles invariant by $\F_0$  around $p_0$ and $\Sigma$ be a transverse 
section 
to $\F_0$ at some point of $\delta_0$. 
Let also $(\C^\mu,\psi)$ be an affine 
chart of $\ffc{d}$ with $\psi(\F_0)=0$. We use also $\psi$ for the points in $\C^\mu$. For instance we denote by $\F_\psi$ the foliation
associated to $\psi\in\C^\mu$ by this affine chart. 

The  holonomy of $\F_0$ along $\delta_0$
in $\Sigma$ is identity. Now considering a $\psi$ near $0$, we have the 
holonomy $h_\psi$ of $\F_\psi$ along $\delta_0$ in $\Sigma$ which is 
called the
deformed holonomy. We write the Taylor expansion
\[
h_\psi(t)-t=\sum_{i=0}^{\infty} a_i(\psi)t^i
\]
The ideal generated by $a_i(\psi), 0\leq i$ is called the Bautin ideal 
of $\delta_0$ in the deformation space $\ffc{d}$. If $\psi\in Zero(I)$ 
then the holonomy of $\F_\psi$ along $\delta_0$ is identity. 
Using Hartogs extension theorem (see also \cite{ily1}), one can see that
the singularity $p_\psi$ near $p_0$ is center 
 and so $\F_\psi\in \wc{d}$. We conclude that $zero(I)\subset\wc{d}$.  

The center $p_0$ of $\F_0$ is called stable if for any deformation $\F_\tau,\tau\in
(\C^k,0)$ of 
$\F_0$ inside $\wc{d}$, the deformed singularity $p_\tau$ is again a center.
Let $\F_0\in\wc{d}$ and $\F_\tau$ be a deformation of $\F_0$ inside $\wc{d}$.
Since each $\F_\tau$ has at least one center, there is a sequence $p_{\tau_i}$ of centers converging to a singularity of $\F_0$. We conclude that the deformed holonomy along the vanishing cycles around $p_0$ is identity and 
$p_\tau$ is center 
for all $\tau$. From this argument we conclude that every $\F_0$ with
$(\wc{d},\F_0)$ irreducible has at least one stable center.
In particular generic points of irreducible components have stable centers.
It is an interesting problem to show that a generic point of $\L(d_1,\ldots,d_s)$ has $d^2-\sum_{i<j} d_id_j$ stable centers.
For the stable center $p_0$ we have $zero(I)=(\wc{d},\F_0)$, where $I$ is the Bautin ideal associated to a vanishing cycle around $p_0$ and the
deformation space $\ffc{d}$.

 Now let $X$ be an irreducible component of $\wc{d}$, 
$\F\in X-sing(\wc{d})$ be a real 
foliation, i.e. its equation has real coefficients, $p$ be a real center
singularity
and $\delta_t,t\in (\R,0)$ be a family of real vanishing cycles around $p$.
The cyclicity of $\delta_0$ in a deformation of $\F$ inside $\ffc{d}$ 
is greater
than $codim_{\ffc{d}}(X)-1$. 
Roughly speaking, the cyclicity of $\delta_0$ is the maximum number of
limit cycles appearing near $\delta_0$ after a deformation of $\F$ in
$\ffc{d}$. The proof of this fact and the exact definition of 
cyclicity
can be found in \cite{rou}. $codim_{\ffc{d}}(\L(d+1))-1=
\frac{(d+2)(d-1)}{2}$ and this is the number obtained by 
Yu. Ilyashenko in \cite{ily}. Now let $X=\L(d_1,\ldots,d_s)$ and $\F=\F(f\sum_{i=1}^s\lambda_i\frac{df_i}{f_i})\in \L(d_1,\ldots,d_s)-sing(\wc{d})$.  Suppose that $\lambda_i$'s and the coefficients of $f_i$'s are real
numbers and $\F$ has a (real) center singularity at $0\in\R^2$. We conclude that
\begin{coro}
Suppose that $s>1$. The cyclicity of $\delta_0$ in a deformation of $\F$ in $\ffc{d}$ 
is not less than 
\[
(d+1)(d+2)-\sum_{i=1}^s (\frac{(d_i+1)(d_i+2)}{2})-1
\]
\end{coro}
Note that the  above lower bound reaches to its maximum when 
$d_1=d_2=\ldots=d_s=1, s=d+1$. In this case the cyclicity of 
$\delta_0$ is not less than $d^2-2$.
Until the time of writing this paper, the best upper bound for the
cyclicity of a vanishing cycle of a  Hamiltonian equation is  the P. 
Mardesic's result $\frac{d^4+d^2-2}{2}$ in \cite{mar}.
 Results for the 
cyclicity of period annulus are obtained by many authors, 
the most complete concerns the Hamiltonian case with $d=2$ 
(see \cite{gav1} and references given there).

We can state center conditions 
for an arbitrary algebraically closed field $k$ instead of $\C$.
The notations in the introduction can be developed for $k$ except 
the center singularity. Suppose that the origin 
$0=(0,0)\in k^2$ is a reduced singularity of $\F(\omega)\in\ffc{d}$. 
It is called a center singularity of $\omega$ if there is a formal power series
$f=xy+f_3+f_4+\cdots+f_n+\cdots$,
where $f_n$ is a homogeneous polynomial of degree $n$ and  with coefficients 
in $k$, such that $df\wedge\omega=0$. 
A singularity $p$ of $\omega$ is called  center if the origin is a center 
singularity of $i^*\omega$, where $i$ is the linear transformation $a\rightarrow a+p$ in $k^2$. Our definition is complete.
Now let $k=\C$ and the origin is a center singularity of a 1-form $\omega$. 
By  theorem A in \cite{mat} the existence of the formal
series $f$ implies the existence of a convergent one, namely $g$. 
Using the  complex Morse theorem we find a coordinates system $(\tilde{x},\tilde{y})$
around the origin such that $g=\tilde{x}^2+\tilde{y}^2$. So our definition
of a center singularity coincides with the definition in  the 
introduction. Now the proof of the fact that $\wc{d}$ is an algebraic subset of 
$\ffc{d}$ is a slight modification of the proof in \cite{mov3}.

\smallskip
\leftline{Hossein Movasati}
\leftline{Institute for studies in
theoretical Physics and Mathematics, IPM} \leftline{School of
Mathematics} \leftline{P.O.Box : 19395-5746}
\leftline{Tehran-Iran}
\leftline{Email: movasati@ipm.ir, hossein.movasati@tu-clausthal.de} 
\end{document}